\documentclass[12pt]{article}
\usepackage{amsmath,amsfonts,amssymb,amsthm,amscd}
\usepackage{a4wide}

\title{\it {Unlikely intersections in   semi-abelian surfaces}}
\author{D. Bertrand
 and H. Schmidt 
\footnote{ { Authors' addresses, AMS Classification, Acknowledgements : see end of the paper.}}}

{Proposition}
{Definition}
{Lemma}
{Corollary}
\newcommand{\C}{\mathbb C}
\newcommand{\R}{\mathbb R}
\newcommand{\Q}{\mathbb Q}
\newcommand{\Z}{\mathbb Z}

\newcommand{\N}{\mathbb N}

\newcommand{\G}{\mathbb G}

\def\mapright#1{\smash{\mathop{\longrightarrow}\limits^{#1}}}

\newcommand{\QQ}{{\mathbb Q}^{alg}}

\date{September 2018 (ArXiv 1803.04835, version 2) }

\begin{document}

\maketitle
 
 \noindent
{\it Abstract}.\,- We consider a   family, depending on a parameter, of  multiplicative extensions of an elliptic curve with complex multiplications. They form a 3-dimensional variety $G$ which admits a dense set of special curves, known as Ribet curves,  which strictly contains the torsion curves.   We show that an irreducible   curve $W$  in $G$ meets this set Zariski-densely only if $W$ lies in a fiber of the family or is a  translate of a Ribet curve by a multiplicative section. 
We further deduce from this result a proof of the Zilber-Pink conjecture (over number fields) for the mixed Shimura variety attached to the threefold $G$, when the parameter space is the universal one. 
 
  \tableofcontents
 
 \vfil \eject
 
\section{Introduction}

\subsection{Statement of the results and plan of the proof}
Let $E_0/{\Q^{alg}}$ be an elliptic curve  
with  complex multiplications. On any  extension  $G_0$  of $E_0$ by $\G_m$ defined over $\Q^{alg}$, there exists a particular subgroup $\Gamma_0$ of $G_0({\Q^{alg}})$, whose elements are called Ribet points. We refer to \S 1.2 below for their precise definition, but point out right now that $\Gamma_0$ contains the  torsion subgroup $G_0^{tor}$ of $G_0(\Q^{alg})$. In fact $\Gamma_0 = G_0^{tor}$ if the extension $G_0$ is isosplit, while $\Gamma_0$ has rank 1 otherwise.

\smallskip
Let further $X/\Q^{alg}$ be a smooth irreducible algebraic curve and let $G/X$ be an $X$-extension   of $E_{0/X}$ by $\G_{m/X}$. Let $q$ be the section   of $\hat E_{0/X} \rightarrow  X$ representing the isomorphism class of the extension $G/X$. We identify $q$ with its image in $E_0(X)$  under the standard polarization $  \hat{E_0} \simeq E_0$, and write $G \simeq G_q$ . Given a section $s$ of $G/X$, we denote by $p = \pi \circ s \in E_0(X)$  its composition with the projection $\pi : G \rightarrow E_0 \times X$. 

\smallskip
Let $\delta \neq 0$ be a  purely imaginary complex multiplication  of $E_0$, and  let $\xi \in X(\Q^{alg})$. A  first property of Ribet points is that if $s(\xi)$ is a Ribet point of its fiber $G_\xi \simeq G_{q(\xi)}$, then its projection $p(\xi)$ to $E_0$ and the point $\delta q(\xi)$ are linearly dependent over $\Z$. Usually, this condition 
alone will be satisfied by infinitely many $\xi$'s. But asking that $s(\xi)$ be a Ribet point   in the fiber of $G_\xi \rightarrow E_0$ above $p(\xi)$ brings a second condition, unlikely to be satisfied infinitely often. And indeed,  we prove in this paper :  
  
 \medskip
\noindent
{\bf Theorem 1} :  {\it  let $G \simeq G_q$ be a non constant (hence non isosplit) extension of $E_{0/X}$ by $\G_{m/X}$, and let $s$ be a section of $G \rightarrow X$, all defined over $\Q^{alg}$. Assume that the set 
 \smallskip

\centerline{ $\mathbf{\Xi} =  \mathbf{\Xi}_s :=  \{\xi \in X(\Q^{alg})$,  $s(\xi)$ is a Ribet point of its fiber $G_\xi \simeq G_{q(\xi)}$\} }

\smallskip
\noindent
 is infinite.  Then,  the sections $p $ and $q$   are  linearly dependent over $End(E_0)$. }

\medskip
 Referring again to \S 1.2  for the definition of the Ribet sections of $G/X$ (which in view of the hypothesis on $G$, also form a  group $\Gamma$ of rank 1, containing the torsion sections), we deduce the following (actually equivalent) version of Theorem 1 :

\medskip
\noindent
{\bf Theorem 2} : {\it assume that the hypotheses of Theorem 1 are satisfied. Then, there exists a non constant or trivial section $s'$ in $\G_m(X)$  such that $s - s'$ is a Ribet section of $G/X$.}

\medskip
The conclusion of Theorem 2 is best possible. Indeed, let $s'$ be such a section in $\G_m(X)$ and let $s''$ be a Ribet section. Then, $s''(\xi)$ is a Ribet point of $G_\xi$ for any $\xi \in X$, while  $s'(\xi)$ lies in $\G_m^{tor}$ infinitely often. The set $\mathbf{\Xi}_s$ attached to $s = s' + s''$ is therefore infinite.

\bigskip
As a corollary to Theorem 2, we consider in an Appendix the case when the curve  $X =  \hat E_0 \simeq Ext(E_0, \G_m)$ is the parameter space of the universal extension ${\cal P}_0$ of $E_0$ by $\G_m$. This extension, which identifies with the Poincar\'e bi-extension of $E_0 \times \hat E_0$ by $\G_m$,  is naturally endowed with the structure of a mixed Shimura variety, for which we prove :

 \medskip
\noindent
{\bf Theorem 3} :  {\it  let $W/\QQ$ be an irreducible algebraic curve in ${\cal P}_0$. Assume that $W$ contains infinitely many points lying on special curves of the mixed Shimura variety  ${\cal P}_0$. Then,  $W$ is  contained in a special surface of ${\cal P}_0$.}

\medskip
Combined with Gao's work on Andr\'e-Oort, this readily implies that the MSV ${\cal P}_0$ satisfies the Zilber-Pink conjecture over $\QQ$ (Conjecture 1.3 of \cite{P}, see \S 5 below).

\bigskip
The proof of Theorem 1 will distinguish three cases. In the first one, we establish the following weaker  version, where the   conclusion is replaced by a ``weakly special" one. Denote by $E_0(\Q^{alg}) \subset E_0(X)$ the group of constant sections of $E_{0/X}$.   

\medskip
\noindent
{\bf Theorem 1.w} : {\it same hypotheses as in Theorem 1. Then,  the sections $p $ and $q$   are  linearly dependent over $End(E_0)$ modulo $E_0(\Q^{alg})$. }

\medskip
The proof of Theorem 1.w   (see \S 2) follows  the  $o$-minimal strategy of Pila-Zannier  and Masser-Zannier, starting with the observation  that if its conclusion does not hold, then the points $\xi$ of $\mathbf{\Xi}$ have bounded height.    

\medskip

In the remaining cases, we suppose that $p $ and $q$   are  linearly dependent over $End(E_0)$ modulo $E_0(\Q^{alg})$. In the second one  (see \S 3), we assume that they are linearly dependent over $\Z$ modulo $E_0(\Q^{alg})$, but that $p$ is not (i.e. $p$ is not constant). Here again, we use the $o$-minimal strategy, but a new argument  is required to check bounded height.

\medskip

In the  last case (see \S 4), we reduce  a weaky special relation over $End(E_0)$ to one over $\Z$, and therefore to a constant section $p$. We finally show that $p$ must be torsion, thanks to a duality argument which turns the problem into a special case of  the Mordell-Lang theorem  (recalled in \S 1.3.(v) below) for a constant  semi-abelian variety attached not to $q$, but to $p$.

\subsection{Ribet sections and points}
 Let  ${\cal X}/\Q^{alg}$ be a smooth irreducible variety, let $A$ be an abelian scheme over  $\cal X$, let $q \in \hat A({\cal X})$ be a section of the dual abelian scheme $\hat A/{\cal X} \simeq Ext_{\cal X}(A, \G_m)$, and let $G = G_q$ be the corresponding ${\cal X}$-extension of $A$ by $\G_{m/{\cal X}}$. We point out that $G_q$ is an isosplit extension  if and only if $q$ is a torsion section. When $A/{\cal X}$ is a constant group scheme, $G_q$ is a constant group scheme  if and only if $q$ is a constant section (for instance a torsion one).

 \smallskip
 
 Let ${\cal P}$ be the Poincar\'e bi-extension of  $A \times_{\cal X} \hat A$ by $\G_m$. For any $\varphi \in Hom_{\cal X}(\hat A, A)$, with transpose $\hat \varphi$, there is 
a canonical isomorphism $ \sigma_{\varphi, q} : {\cal P}((\varphi  - \hat \varphi) (q) , q)  \simeq \G_{m/{\cal X}}$ of $\G_m$-torsors over $\cal X$ (see \cite{ACL} Prop. 6.3, \cite{PcPbGj}). We define the {\it basic Ribet section}  associated to $\varphi$  as the section  
$s_{\varphi,q}  = \sigma_{\varphi, q}^*(1_{\cal X})$ 
 of the semi-abelian scheme $G = G_q = (id_A, q)^*{\cal P} = {\cal P}_{| A \times q}$ over $\cal X$. We say ``point" instead of ``section"  if ${\cal X}$ is a point, and drop the index $q$ when the context is clear.  
 
 \smallskip
 The Ribet section $s_\varphi \in G({\cal X})$ depends additively on $\varphi$ (cf. \cite{JR}, Prop. 4.2; \cite{PcPbGj}), and in fact only on  $\varphi - \hat \varphi$.  Its projection under $\pi : G \rightarrow  A$ is the section 
 $$p_\varphi :=  \pi \circ q_\varphi = (\varphi - \hat \varphi) \circ q \in A({\cal X}).$$
 So, when $\varphi$ varies, the basic Ribet sections form a finitely generated  subgroup of $G({\cal X})$, of rank $r_q$ at most equal to the rank of the $\Z$-module ${\cal E} = \{\varphi  - \hat \varphi,  \varphi \in Hom_{\cal X}(\hat A, A)\}$, and equal to it when $q$ is sufficiently general. On the other hand, $r_q = 0$ if $q$ is a torsion section. Indeed,  although their dependence in $q$ is {\it not linear}, the Ribet sections $s_\varphi$  satisfy the following ``lifting property" (for $(i) \Rightarrow (ii)$, see \cite{B-E}, \cite{BMPZ}  in the case of points, and \cite{PcPbGj} in general) :

\medskip
\noindent
{\bf Lemma 1} : {\it let $\varphi \in Hom_{\cal X}(\hat A, A)$,  let $q \in \hat A({\cal X})$ and consider the conditions

(i)  $q$  is a torsion section

(ii)   $s_\varphi$ is a torsion section

(iii) $p_\varphi$ is a torsion section.

\noindent
Then, $(i) \Rightarrow (ii) \Rightarrow (iii)$, and if $\varphi - \hat  \varphi$ is an isogeny, the three conditions are equivalent. }  

\medskip
More generally, let $s$ be a local section of $G \rightarrow {\cal X}$ (for the \'etale topology). We say that $s$ is a {\it Ribet section}   of  $G/{\cal X}$ if there exists a positive integer $n$ satisfying\;: $n.s = s_\varphi$ for some $\varphi$, with multiplication by $n$ in the sense of the group scheme $G/{\cal X}$. The projection $p$ of $s$ to $A$ satisfies : $n p = (\varphi - \hat \varphi)\circ q$.  All (local) torsion sections of $G/{\cal X}$ now appear as such Ribet sections, and Lemma  1  extends to this more general setting. Viewed as points above the generic point $\eta$ of ${\cal X}$, with $K = \Q^{alg}({\cal X}_\eta)$, the Ribet sections form a subgroup $\Gamma$  of the group $G_\eta(K^{alg})$, of same rank $r_q$ as above. 
 
\medskip
The construction of Ribet sections commutes with  any base change. For instance,  given a  basic Ribet section $s_{\phi, q}$ of $G/{\cal X}$, and a point $\xi$  in ${\cal X}(\Q^{alg})$,  $s_{\varphi, q}(\xi) = s_{\varphi_\xi, q(\xi)}$ is the basic Ribet point of the fiber $G_\xi$ attached to the specialization $\varphi_\xi$ of $\varphi$ at $\xi$.   Conversely, let $s^\xi$ be a Ribet point of $G_\xi(\Q^{alg})$. By definition, there exist $n_\xi \in \Z_{>0}$ and $\varphi_\xi \in Hom(\hat A_\xi, A_\xi)$ such that $n_\xi s^\xi = s_{\varphi_\xi, q(\xi)}$. {\it Assume} further {\it that  $\varphi_\xi$ extends to an element $\varphi \in Hom(\hat A, A)$} (which occurs automatically if $A/{\cal X}$ is a constant abelian scheme as in \S 1.1). Then,  $s_{\varphi_\xi, q(\xi)} = s_{\varphi, q}(\xi)$, and  there exists  a local section $s$ of $G / {\cal X}$  such that  $n_\xi.s = s_{\varphi}$, whose image in $G$ contains  $s^\xi$. So, the Ribet point $s^\xi$  extends locally to a Ribet section of $G/{\cal X}$.   

\bigskip 

Let us now return to the situation of \S 1.1,  where $A = E_0 \times {\cal X}$, for a CM elliptic curve $E_0$, and ${\cal X}$ is either the curve $X$ or a point $\xi$ on $X$. Then, the $\Z$-module $\cal E$ above identifies with
$${\cal E} = \{\varphi  - \overline \varphi,  \varphi \in   End(E_0) \} = \Z \delta,$$
where $\delta = \alpha - \overline \alpha \neq 0 $ is a purely imaginary quadratic number, which will be fixed from now on. Consequently, for any $q \in E_0(X)$,  the group of basic Ribet sections of $ G = G_q$ is cyclic, generated by the section
$$s^R := s_{\alpha, q} \in G(X) \; , \;  ~�{\rm {with}}�~ p^R := \pi \circ s^R = \delta q \in E_0(X).$$
Viewed at the generic point $\eta$ of $X$, {\it the Ribet sections of $G/X$ then form the divisible hull $\Gamma$  of the group $\Z. s^R(\eta)$ in $G_\eta(K^{alg})$}. Furthermore, 
for any $\xi \in X(\Q^{alg})$, the value $s^R(\xi) = s_{\alpha, q(\xi)}$ of $s^R$ at $\xi$ generates the group of basic Ribet sections of $G_\xi = G_{q(\xi)}$, and the Ribet points of $G_\xi$ form the divisible hull  
$$\Gamma_\xi=\{ s^\xi \in G_\xi(\Q^{alg}), \exists (n, m) \in \Z^2, n \neq 0 , ns^{\xi} = m s^R(\xi)\}  \supset G_\xi^{tor}$$
of $\Z.s^R(\xi)$ in  $G_\xi(\Q^{alg})$. 

\medskip
Under the assumptions of \S 1.1, the section $q$ is not constant, hence not torsion, while $\delta$ is an isogeny, so $s^R$ is not torsion by Lemma 1, and the rank $r_{q}$ of $\Gamma$ is  equal to 1. On the other hand, by Lemma 1 (now at the level of points), given a point $\xi \in X(\Q^{alg})$, 
$$q(\xi) \in E_0^{tor}  \Leftrightarrow s^R(\xi) \in G_\xi^{tor} \Leftrightarrow \Gamma_\xi = G_\xi^{tor},$$
and this  occurs for infinitely many $\xi$'s since $q$ is not constant 
 (cf. \cite{B-E}, Thm. 1). Otherwise, $\Gamma_\xi$  has rank 1, but for $s(\xi) \in \Gamma_\xi$, we still have : $s(\xi) \in G_\xi^{tor} \Leftrightarrow p(\xi) \in E_0^{tor}$. 

\medskip
In view of these descriptions of the groups $\Gamma$ and $\Gamma_\xi$, our work can be interpreted as a  
particular case of the study of unlikely intersections within an isogeny class (cf. \cite{Gao}), or of a relative version of the Mordell-Lang problem (compare with \S 1.3.(v) below).

 \subsection{The context}
We here put the results of \S 1.1 in perspective with other statements of unlikely intersections.  Two sets 
$$   \mathbf{\Xi}^{tor}   \subset   \mathbf{\Xi}  \subset   \mathbf{\Xi}^{\ell d}$$
related to the section $s \in G(X)$  naturally appear in the process.

\medskip

(i) Theorem 1 gives a positive answer to the ``Question 2" raised in \cite{NDJFL}, \S 5, while  a positive answer to its ``Question 1" was recently obtained by Barroero  \cite{Ba}. However, the applications to Pink's conjecture given in \cite{NDJFL} require  clarification, because of their ambiguous use of  Hecke orbits.  
We by-pass this problem for the mixed Shimura variety  ${\cal P}_0$ studied in the Appendix, by describing all its possible special curves. Theorem 3 will then follow from Theorem 2,  along the method of \cite{NDJFL}.  

\medskip

(ii) Contrary to the convention of \cite{BMPZ},  the torsion points are here viewed as particular cases of Ribet points. Therefore,  Theorem 2  implies the restriction to the case of our semi-abelian scheme $G/X$ of the main theorem of \cite{BMPZ}, which concerns the subset 

\smallskip
\centerline{ $\mathbf{\Xi}^{tor} = \mathbf{\Xi}^{tor}_s :=  \{\xi \in X(\Q^{alg})$,  $s(\xi)$ is a torsion point of its fiber $G_\xi$ \}}

\smallskip
\noindent
of $\mathbf{\Xi}$, and asserts :

\medskip
\noindent
{\bf Lemma 2} :  {\it let $G/X$ and $s$ be as in Theorem 1, and assume that  the set  $\mathbf{\Xi}^{tor}$  is infinite. Then $s$ is  a Ribet section or  a torsion translate of a non constant section in $\G_m(X)$.}

\medskip
For $\xi \in \mathbf{\Xi}^{tor}$,  $p(\xi)$ too is torsion, so (by  the Manin-Mumford theorem \cite{H} for the image of $(p, s')$ in $E_0 \times \G_m$), the conclusion of Theorem 2 can be sharpened to  the same statement. 

\medskip
We pointed out at the end of \S 1.2 that the  set $ \mathbf{\Xi}^{tor}_{s^R}$ attached to $s^R$ (hence more generally to any Ribet section) is infinite. Therefore, Lemma 2 too is best possible.

\medskip

(iii)  In relation with the two sections $s, s^R$ of $G/X$, consider the set

\smallskip

\centerline{ $\mathbf{\Xi}^{\ell d} = \mathbf{\Xi}^{\ell d}_{s, s^R} :=  \{\xi \in X(\Q^{alg})$,  $s(\xi)$ and $s^R(\xi)$ are linearly dependent over $\Z$\}   . }

\smallskip
\noindent
For $\xi$ in this set, either $s(\xi)$ lies in the divisible hull $\Gamma_\xi$ of $\Z. s^R(\xi)$, or $s^R(\xi)$ is a torsion point. So  $\mathbf{\Xi}^{\ell d}$ is the  union of  $\mathbf {\Xi}$ and  $ \mathbf{\Xi}^{tor}_{s^R}$ and in particular, is  always infinite.  More generally, given two sections $s, s'$  in $G(X)$, the similarly defined set $\mathbf{\Xi}^{\ell d}_{s,s'}$ will be infinite as soon as the group generated by $s$ and $s'$ in $G(X)$ contains a non-torsion Ribet section. So, in contrast with the case of abelian schemes (see \cite{MZ}, \cite{B-C}), the subgroup schemes of $G \times_X G$ do not suffice to control  the finiteness of $\mathbf{\Xi}^{\ell d}_{s,s'}$ : as in \cite{PcPbGj}, the special subvarieties of the corresponding mixed Shimura variety  should also be taken into account.
 
\medskip

(iv) Consider the curve $W = s(X)$ in $G$ and define a Ribet curve as  the image in $G$ of  a Ribet section. Theorem 2  then says that $W$ is the translate of a Ribet curve by a section in $\G_m(X)$. Since any curve $W$ in $G$ dominating $X$ can be viewed as the image of a section after a base extension, while any Ribet point of a fiber $G_\xi$  locally extends to a Ribet section, this justifies the
 last but one sentence of the abstract.

\medskip

(v) Assume that contrary to the hypothesis of Theorem 1 ,  $G = G_0 \times X$ for some constant semi-abelian surface $G_0/\Q^{alg}$,  and that $s$ is not constant. Then, the projection $W_0$ of $W = s(X)$ to $G_0$ is  a curve, which contains infinitely many points of the group $\Gamma_0$ of Ribet points of $G_0$. Since $\Gamma_0$ has finite rank (at most 1),  the solution by Vojta and McQuillan \cite{McQ} of the {\it Mordell-Lang} conjecture  for semi-abelian varieties implies that $s$ factors through a translate by a Ribet point of  a strict connected algebraic subgroup of $G_0$. If the section $q$, here constant, is not torsion, the only such one is $\G_m$.  So the conclusions of the theorems still hold  true in this case.

\medskip

(vi) Same as in (v), but assume furthermore that $q$ is a torsion section, say the trivial one, so $G_0 \simeq \G_m \times E_0$.  Then, $s = (s', p)$ for some section $s' \in \G_m(X)$, while the group $\Gamma_0$ of Ribet points of $G_0$   coincides with   $G_0^{tor}$.  By Manin-Mumford,  $\mathbf{\Xi} = \mathbf{\Xi}^{tor}$ is then infinite if and only if $s'$ is a torsion section, or $p$ is a torsion section. 

\medskip

(vii) In this paper, we do not touch on the question of replacing $\Q^{alg}$ by $\C$, or of applying Theorem 2  to generalized Pell equations as in \cite{MZ}, \cite{B-C}. Nor do we study  how effective our results can be made.  Note that Lemma 2 above is made effective in the  ongoing work  \cite{JS2}. Due to the use of Pfaffian methods, in particular \cite{JT} and \cite{JS1}, the bounds for the counting problem in \cite{JS2} are uniform and effective.  

\medskip
We take opportunity of these comments to show that 

\smallskip
\noindent
{\it Theorem 1 $\Leftrightarrow$ Theorem 2} : Theorem 2 clearly implies Theorem 1. Indeed,   the sections $s$ and  $s'' = s - s'$ have the same projection $p$ to $E_0$. Since $s''$ is a Ribet section,  $p$  and $\delta q$ are linearly dependent over $\Z$, so $p$ and $q$ are linearly dependent over $End(E_0)$. 

\smallskip

Conversely, assume that the hypotheses and the conclusion of Theorem 1 hold true, and let $np - \rho \, q = 0$ be a non-trivial relation with $n \in \Z, \rho \in End(E_0)$ not both $0$ (equivalently, $n \neq 0$ since $q$ is not a torsion section). Without loss of generality, we can assume that $\mathbf{\Xi}^{tor}$ is finite, otherwise Lemma 2 readily implies the conclusion of Theorem 2.  For any $\xi \in   \mathbf {\Xi}$, $\delta q(\xi)$ and the the projection $p(\xi)$ of the Ribet point $s(\xi)$ are linearly dependent  over $\Z$,  so there exist  $n_\xi , m_\xi \in \Z$,  not both zero, such that $n_\xi p(\xi) - m_\xi \delta  q(\xi) = 0$, while the generic relation implies : $np(\xi) - \rho \, q(\xi) = 0$. If these two relations are  linearly independent over $End(E_0)$, then $q(\xi)$, hence $s^R(\xi)$, hence $s(\xi)$, are torsion points and $\xi$ lies in $\mathbf{\Xi}^{tor}$. So, for infinitely many, hence at least one,  $\xi$,  these two relations must be linearly dependent  over $End(E_0)$, and in fact over $\Z$, since $n $ does not vanish. This implies that $\rho$ is a rational multiple of $\delta$, and by their very construction, this in turn  implies the existence of a Ribet section $s''$ projecting to $p$. So, $s' = s - s''$ factors through $\G_m$. Finally, if $s'$ is a  constant section, it must be a torsion one since $s'(\xi)$   is a Ribet point  of $G_\xi$ projecting to $0$ for one (any) $\xi \in \mathbf{\Xi}$. In this case, $s$ itself is a Ribet section, and otherwise $s'$ is not constant, so the conclusion of Theorem 2 holds in all cases. \qed

\section{Proof of Theorem 1.w}
Recall the hypotheses of Theorem 1.w,  as well as the notation $s^R, \Gamma_\xi,...$ of \S 1.2.  So, $q \in E_0(X)$ is not constant, $s$ is a section of $G= G_q \rightarrow X$ projecting to the section $\pi \circ s = p \in E_0(X)$, and the set $\mathbf \Xi = \{�\xi \in X(\Q^{alg}), s(\xi) \in \Gamma_\xi \}$, concretely described as  
$$ \mathbf \Xi = \{\xi \in X(\Q^{alg}), \exists (n, m) \in \Z^2,  n \neq 0, ns(\xi) - m s^R(\xi) = 0\}.$$
is infinite.   {\it We assume that the sections $p$ and $q$ are linearly independent over $End(E_0)$ modulo $E_0(\Q^{alg})$}, and search for a contradiction. 

\smallskip

 We fix  a number field $k$ over which $X$ and $G$, hence the sections $q$ and $s^R$, as well as the section $s$, hence $p$, and the isogeny $\delta$, are defined. We recall that the basic Ribet section $s^R$ projects to $E_0$ on the section $p^R =   \delta q $.

  \subsection{The o-minimal strategy}
The proof of Theorem 1.w will be done in 5 steps. The 3rd one is developed in \S 2.2. By a ``constant" $c , \gamma $, we mean a positive real number  which depends only on the data $X, E_0, q, s$ and the number field $k$. The constants $C $ may depend on further data introduced in the proof. 

\smallskip
We point out that any finite set of points can without loss of generality be withdrawn from the curve $X$. To ease a technical point in the 3rd step, we will for instance require that the sections $p, q$ and $p+q \in E_0(X)$ never vanish on $X$. The complement is a finite set since  $q$ is not constant, $p$ can be assumed to be so (constant $p$'s  are treated by a direct method in \S 4.2), and if $p+q$ is constant, we can make it non constant by replacing $s$ by $2s$, so $p$ by $2p$, without modifying the content of the theorems.    
 
\medskip
\noindent
{\bf 2.1.1. Bounded heights of points}

\medskip
Let $h$ denote a height on $X(\Q^{alg})$ attached to a divisor of degree 1 on the completed curve. Consider the set

\medskip
\centerline{ $\mathbf{\Xi}^{\Z\ell d}_{p, \delta q} = \{\xi \in X(\Q^{alg}), p(\xi)$ and $\delta q(\xi)$ are linearly dependent over $\Z$\}.} 

\medskip
\noindent
Since the projection $p(\xi) = \pi \circ s(\xi)$ of a Ribet point $s(\xi)$ lies in the divisible  hull of the group $\Z. \delta q(\xi)$  in $E_0(\Q^{alg})$, this set contains ${\mathbf \Xi}$.

\medskip 
\noindent
{\bf Lemma 3}�: {\it let $p, q \in E_0(X)$ be linearly independent over $End(E_0)$ modulo $E_0(\Q^{alg})$. There exists a constant $c_0$ such that  $h(\xi) \leq c_0$ for any $\xi \in \mathbf{\Xi}^{\Z\ell d}_{p, \delta q} $, and in particular, for any $\xi \in \mathbf{\Xi}$.}

\smallskip
\noindent
{\it Proof}. - In view of the hypothesis on $p,   q$, bounded height on $\mathbf{\Xi}^{\Z\ell d}_{p, \delta q}$ follows directly  from \cite{V}, Theorem 4 (and one can even replace $\Z$ by $End(E_0)$ in the definition of $\mathbf{\Xi}^{\Z\ell d}_{p, \delta q}$). Alternatively, one can appeal to Silverman's specialization theorem \cite{S}. 
   \qed

\smallskip
To get the searched for contradiction, it remains to show that the degrees 
$$d_\xi = [k(\xi):\Q]$$ 
too are bounded from above on the set  $\mathbf \Xi$. 

\medskip
\noindent
{\bf 2.1.2. Heights of relations bounded by degrees}

\medskip 
\noindent
{\bf Lemma 4}�: {\it there exist  two constants $c \, , \gamma$ such that  for any point $\xi \in {\mathbf{\Xi}}$, there exist  two integers $n \neq 0, m$ of absolute values $\leq c d_\xi^\gamma$ such that $n s(\xi) - m s^R(\xi) = 0 $.}

\smallskip
\noindent
{\it Proof}. -  By \cite{BMPZ}, Corollary of \S 3.1, there exists a constant $c'$ such that if $s(\xi)$ is a torsion point of $G_\xi$,   its order $n$ is bounded from above  by $c' d_\xi^4$, so $(n,0)$ satisfies the required condition. 
We can therefore assume that the Ribet point $s(\xi)$, hence $q(\xi)$ by Lemma 1, is not a torsion point. For $\xi \in {\mathbf{\Xi}}$, there exist $a, b \in \Z$, not both $0$,  such that $ap(\xi)  - b \delta q(\xi) =0$, and since $q(\xi) \notin E_0^{tor}$, any such relation will automatically have $a \neq 0$. The points $p(\xi), \delta q(\xi)$ are defined over $k(\xi)$, and have heights  $\leq c_0$. By works of Masser and  David, 
there then exists such a relation with $max (|a|, |b|) \leq c_1 d_\xi^{\gamma_1}$ for some constants $c_1, \gamma_1$. 

\smallskip

By our running hypothesis that  $q(\xi)$ is not  torsion, the set of such relations (trivial one included) is a free $\Z$ module of rank 1, and its generator $(a_0, b_0)$ satisfies the above bound. 

\smallskip
Consider now the non-torsion Ribet point $s(\xi)$ (so,  $s^R(\xi)$ too is non-torsion), and let  $(n_0 \neq 0, m_0) \in \Z^2$ be a generator of the group of relations  $n s(\xi) -  m s^R(\xi) = 0$, which is again free of rank 1.  Projecting to $E_0$, we then have  $n_0p(\xi) - m_0 \delta q(\xi) = 0$.  So, there exists $d \in \N$ such that $(n_0, m_0) = d. (a_0,b_0)$, and  $a_0s(\xi) - bs_0^R(\xi)$ is a  torsion point of $G_{q(\xi)}$, of exact order $d$ since $(n_0, m_0)$ is  minimal. Since it projects to $0$ on  $E_0$, it is  actually a $d$-th root of unity $\zeta_d$ .  Now, both $s(\xi)$ and $s^R(\xi)$ are defined over $k(\xi)$ (since  $s$ and $s^R$ are global sections of $G \rightarrow X$), so $\zeta_d$ too lies in $k(\xi)$. Since $\zeta_d$ has order $d$,  this implies that $d \leq c_2 d_\xi^{\gamma_2}$, say with $\gamma_2 = 2$. 

\smallskip
In conclusion, for any $\xi \in {\mathbf{\Xi}}$, there is a linear relation $ns(\xi) - ms^R(\xi) =0$, with $(n,m) \in\Z^2, n \neq 0$ and $max(|n|, |m|) \leq c d_\xi^\gamma$ for some constants $c$ and $\gamma = \gamma_1 + \gamma_2$.  \qed

\medskip
\noindent
{\bf 2.1.3. Counting relations of bounded height} 

\medskip
\noindent
In this step and the next one, we extend the scalars from $\Q^{alg}$ to $\C$, but still write  $X, K = \C(X)$, etc, instead of $X_\C, K \otimes \C, ...$. We sometimes indicate by the exponent $^{\it an}$   the analytic object attached to an algebraic one over $\C$. 

\smallskip
 We now follow the usual procedure of studying the lifts to a universal covering  of the relations considered in Lemma 4,  and bounding their number via (generalizations of) the Pila-Wilkie theorem for a relevant $o$-minimal structure.  There are several ways to implement this method. For instance, we can 
 
 \medskip
\noindent
(A) choose a fundamental domain ${\cal F}$ for the uniformization map ${\rm unif} :  \tilde {G} \simeq  \C \rtimes (\C  \times \tilde X)  \rightarrow G^{an}$,  and count the relations in $\tilde{G}$ when  the transcendence degree  over $\C$ of the field of definition of 
${\rm ({unif}_{|{\cal F}})^{-1}} \circ s$ is large enough. Here,  ${\cal F}$ is unbounded, but by work of Peterzil and Starchenko, a convenient choice allows to work in the $o$-minimal structure $\R_{an, \exp}$ ; 
 
\smallskip
\noindent
(B) or fix  a simply connected domain $D \subset X^{an} $,  consider the exponential morphism $\exp_G$,  restricted over  $D$, and count the relations in $(Lie G)/D \simeq (\C \rtimes \C) \times D$ when the transcendence degree over $\C(X)$ of the field of definition of  $\exp_G^{-1} (s_{|D})$ is sufficiently large. Here, $D$ can be compact, and it suffices to work in the $o$-minimal structure $\R_{an}$. 

\smallskip
An advantage of (A) is its impact on effectivity, as alluded to in Comment (vii) of  \S 1.3  (see also Remark 3 of \S 4.3).  
But as in \cite{BMPZ}, \S3.3, we here follow the more elementary approach (B), taking advantage of the computation of transcendence degrees already established in this paper.

\medskip
So, let $(D, \xi_0)$ be a pointed set in $X^{an}$, homeomorphic to a closed disk.  The group scheme $G/X$ defines an analytic family $G^{an}$ of Lie groups over the Riemann surface $X^{an}$. Similarly, its relative Lie algebra $(LieG) /X$ defines an analytic vector bundle $LieG^{an}$ over $X^{an}$, of rank 2. We denote by $\Pi_G$ the  $\Z$-local system  of periods   of $G^{an}/X^{an}$; it   is the kernel of the  exponential exact sequence of analytic sheaves over $X^{an}$ :
$$0~ \mapright{} ~{\Pi}_G~ \mapright{} ~Lie G^{an} ~\mapright{\exp_G}�~ G^{an} ~\mapright{}~ 0 ~,.$$
For any $U_0$ in $Lie(G_{\xi_0}(\C))$   such that $\exp_{G_{\xi_0}}(U_0) = s(\xi_0) \in G_{\xi_0}(\C)$, there exists a unique analytic section $U$ of $Lie(G^{an})/D$ (meaning : over a neighbourhood of $D$), such that 
 $$U(\xi_0) = U_0 ~�{\rm ~�and~}�~\forall \xi \in D, \exp_{G^{an}_\xi}(U(\xi)) = s(\xi). $$
 Since $D$ is fixed, we will
just write $U = \log_G(s)$, although only its class modulo $\Pi_G$ is well defined. Similarly, let $U^R = \log_G (s^R)$ for the Ribet section $s^R$.  By the same process for $E_{0/X}$ (and  the tacit assumption that the logarithms at $\xi_0$ are chosen in a compatible way), the projection $p = \pi \circ s \in E_0(X)$ admits as logarithm
 $\log_{E_0} ( p) := u  = d\pi(U)$; we also set $v = \log_{E_0}(q)$, so $d\pi(U^R) := u^R =   \delta v$. 
 
 \smallskip
 We will use the explicit expressions given in \cite{BMPZ} for $U, U^R$ and $\Pi_G$. These hold on any simply connected domain of $X^{an}$ where $u, v$ and $u+v$ do not assume period values. This is ensured by the hypothesis, made at the beginning of \S 2.1, that $p, q$ and $p+q$ vanish nowhere on $X$.
  
 \medskip
Let $K = \C(X)$ be the field of rational functions of $X$. Since $Lie G$ is a vector bundle over $X$, it makes sense to speak of the field of definition $K(U)$ of $U$ over $K$. 
Similarly, let $F_G = K(\Pi_G)$ be the field of definition of   $\Pi_G$. Notice that the field $F_G(U)$ now depends only on the section $s$. Moreover, for the Ribet section $s^R$, we have : 
  
  \medskip
  \noindent
  {\bf Lemma 5}�: {\it  the field of definition $F^R = K(U^R)$ of any logarithm $U^R$ of $s^R$  coincides with the field of periods $F_G$ of $G$.} 
   
   \medskip
   \noindent
   {\it Proof}. -  The  explicit expressions of  $\Pi_G$ and $U^R$ given in \cite{BMPZ}, \S A.1, show  that both fields coincide with the field $K(v, \zeta(v))$, where $\zeta$ denotes the Weierstrass zeta function of the elliptic curve $E_0$. \qed
   
\medskip
For any real number $T \geq 1$, set $\Z[T] = \{n \in \Z, |n|�\leq T \}$, and consider the subset
$${\mathbf \Xi}[T] := \{�\xi \in X(\Q^{alg}), \exists (n, m) \in (\Z[T])^2, n \neq 0, ns(\xi) - m s^R(\xi) = 0 \}$$
of ${\mathbf \Xi} = {\mathbf \Xi}_s$. We then have :

\medskip
\noindent
{\bf Proposition 1}�: {\it let $D$ be a closed  disk in $X^{an}$. For any $\epsilon > 0$, there exists a  real number $C_\epsilon$\,, depending only on $X, E_0, q, s,  D$ and $\epsilon$\,, such that

\medskip

- {\rm (a)}   either, for any $T \geq 1$, there are at most $C_\epsilon T^\epsilon$ points in $D \cap {\mathbf{\Xi}[T]}$;

\smallskip

- {\rm (b)}�  or the field $F_G(U)$ has transcendence degree at most 1  over the field $F_G$.}

\medskip
\noindent
The proof of Proposition 1 is given in \S 2.2 below, as a corollary of Habegger-Pila's ``semi-rational" count \cite{HP}, Corollary 7.2.

   \bigskip
   
  \noindent
  {\bf 2.1.4. Logarithmic Ax}
  
  \smallskip

Assume that Conclusion (b) of Proposition 1 holds.  Since $u = d\pi(U)$, the field $F_G(U)$ has transcendence degree at most 1 over  $F_G(u)$, and  :

\smallskip

- (b1) either $u$ is algebraic over $F_G = K(v, \zeta(v))$, in which case we know by the Ax-Schanuel theorem on (the universal extension of) the  elliptic curve $E_0$ that $p$ and $q$ are linearly dependent over $End(E_0)$ modulo constants;    

- (b2) or   $U = \log_G(s)$ is algebraic over $F_G(u)$, hence over $K(u, \zeta(u), v, \zeta(v))$, in which case we know by \cite{BMPZ}, Lemma 5.1,  that $s$ is a translate of a Ribet section by a constant  one, i.e. one in $\G_m(\C)$
 since $G$ is not isosplit.  Then, $p = \pi \circ s$ and $q$ are linearly dependent over $End(E_0)$.

\smallskip
In both cases, we get a contradiction to our hypothesis that $p$ and $q$ are linearly independent over $End(E_0)$ modulo $E_0(\Q^{alg})$. So, Conclusion (a) must hold.

\medskip
  \noindent
  {\bf 2.1.5. Conclusion}
  
  \medskip

It follows from Lemma 3 and  a compactness  argument (see \cite{MZ}, Lemma 8.2 and the paragraph after (9.2)) that there exists  a finite set of closed disks $D_i$ in $X^{an}$ and a constant $c'$ such that the following holds : for any $\xi \in {\mathbf \Xi}$,  a positive proportion  $\frac{1}{c'} d_\xi$ of the conjugates of $\xi$ over $k$ lie in one of the $D_i$'s, say $D_1$. Now,  all these conjugates are    still in ${\mathbf{\Xi}}$, since $\sigma(s^R(\xi)) = s^R (\sigma \xi)$ is a Ribet point of $G_{q(\sigma \xi)}$ for $\sigma \in Gal(\Q^{alg}/k)$.  Actually, by Lemma 4, all the  conjugates of $\xi$  over $k$ lie  in ${\mathbf \Xi}[T]$ with   $T = c d_\xi^\gamma$.  Choosing $\epsilon = 1/2\gamma$, we deduce from Conclusion (a) that $D_1 \cap {\mathbf{\Xi}}$  has at most $c''d_\xi^{1/2}$   (and  at least $\frac{1}{c'} d_\xi$\,) elements.   Therefore, $d_\xi$ is bounded from above on ${\mathbf{\Xi}}$, and this concludes the proof of Theorem 1.w.

\subsection{The semi-rational count}  
The proof of Proposition 1 uses Betti coordinates and maps, defined as follows. We recall that $D \subset X^{an}$ is homeomorphic to a closed complex disk. 

 \medskip
 The sections of the local system $\Pi_G$ over $D$ form a $\Z$-module $\Pi_G(D) \subset Lie G^{an}(D)$ of rank 3, with a basis $\{\varpi_0, \varpi_1, \varpi_2\}$ such that $\varpi_0$ generates $\Pi_{\G_m}(D)$, and $\varpi_1, \varpi_2$ project to a basis $\omega_1, \omega_2$ of $\Pi_{E_0}(D)$.   Then, any logarithm $U := \log_G (s)$ of a section $s$ of $G/X$  over the disk $D$ can uniquely be written as  
 $$U = b_0 \varpi_0 + b_1 \varpi_1 + b_2 \varpi_2,$$
  where $b_0, b_1, b_2$ are real analytic functions on $D$, with values in $\C$ for $b_0$, and in $\R$ for $b_1$ and $b_2$. We call $(b_0, b_1, b_2)$ the Betti coordinates of $U$, and define the Betti map attached to $U$ as
 $$ U_B = (b_0; b_1, b_2) : D \rightarrow \C \times \R^2,$$
 Similarly,  we write $U^R_B = (b^R_0; b^R_1, b^R_2) $ for the Betti map attached to $U^R = \log_G(s^R)$, and denote by 
 $\cal S$ the image of the disk $D$ under the map 
$${\cal U}_B := (U_B, U_B^R): D \twoheadrightarrow {\cal S} \subset  \R^4 \times \R^4 = \R^8.$$
We will work in the $o$-minimal structure $\R_{an}$ of globally subanalytic sets. 

\medskip
\noindent
{\bf Lemma 6}� : {\it ${\cal S} = {\cal U}_B (D)$ is a compact 2-dimensional set, definable in the structure $\R_{an}$.}

\smallskip
\noindent
{\it Proof}. - By definition (or by inspection of the formulae in \cite{BMPZ}), the maps $U_B$ and $U^R_B$ extend  to real analytic maps on a neighbourhood of the compact disk $D$. Therefore, ${\cal S} = {\cal U}_B (D)$ is  a compact definable set. Furthermore, the Betti map  $\pi \circ U^R_B := u^R_B = (b^R_1 ,  b^R_2)$  attached to $u^R = \log_{E_0}(p^R)$  is an immersion (since $p^R = \delta q \in E_0(X)$ is not a constant section), so $\cal S$ is indeed a real surface.  \qed

 \medskip
 With this notation in mind, a point $\xi$ of  $D$ lies in $D \cap \mathbf {\Xi}$ if and only if 
 $$ \exists (\nu \neq 0, \mu) \in \Z^2, \exists (\beta_0, \beta_1, \beta_2) \in \Z^3, \nu U(\xi) - \mu U^R(\xi) = \beta_0 \varpi_0(\xi) + \beta_1 \varpi_1(\xi) + \beta_2 \varpi_2(\xi),$$
or alternatively,  in terms of the Betti maps :
$$ \exists (\nu \neq 0, \mu) \in \Z^2, \exists (\beta_0, \beta_1, \beta_2) \in \Z^3, \nu U_B(\xi) - \mu U^R_B(\xi) = (\beta_0; \beta_1, \beta_2) \in \Z \times \Z^2 \subset \C \times \R^2.$$
Remark  that 

- {\it if $|\nu|, |\mu|$ are bounded by some number $T$, then   $|\beta_0|, |\beta_1|, |\beta_2| \leq  C_1 T$} for some constant $C_1$, since $D$ is compact; 

-   {\it given any}  real {\it numbers $\nu \neq 0, \mu, \beta_0, \beta_1, \beta_2$}, there are only {\it finitely many $\xi$'s in $D$ such that $\nu U_B(\xi) - \mu U^R_B(\xi) = (\beta_0; \beta_1, \beta_2)$}. Otherwise,  $\nu u - \mu \delta v $ would be constant on $D$, contradicting the Ax-Schanuel theorem invoked in \S 2.1.4.(b1).

\bigskip 
We can now describe the definable set $\cal Z$ to which Habegger-Pila's semi-rational count \cite{HP} will be applied. On the one hand, we have the affine space $\R^5$  with real coordinates $(\nu, \mu, \beta_0, \beta_1, \beta_2)$; we will indicate by the index $_*$ the complement of the hyperplane $\nu = 0$.  On the other hand, we have the affine space $\C \times \R^2 = \R^4$ and its square $\R^8$, which is the target space of the map ${\cal U}_B$.   We  consider the incidence variety $\cal Z$ in $\R^5 \times \R^8$, with projections $\pi_1$ to $\R^5_* \subset \R^5$ and $\pi_2$ to ${\cal S} =  {\cal U}_B (D) �\subset \R^8$ :
$${\cal Z} = \{\big( (\nu, \mu, \beta_0, \beta_1, \beta_2) ; \big(w := (w_0; w_1, w_2), w^R := (w_0^R; w_1^R, w_2^R)\big)\big) \in \R^5 
\times {\cal S} \subset \R^5 \times \R^8, $$
$$ \qquad \qquad {\rm such ~�that} ~\nu \neq 0 ~{\rm and}  ~ \nu. w -\mu. w^R =  (\beta_0; \beta_1, \beta_2)  \in \R \times \R^2 \subset \C \times \R^2 = \R^4\}$$
By Lemma 6, ${\cal Z}$ is a definable subset of $\R^{13}$. Furthermore, ${\cal U}_B(D \cap \mathbf \Xi) = \pi_2(\pi_1^{-1}(\Z^5_*))$.

\medskip

Let $\epsilon \in \R_{>0}$. Given $T \geq 1$,  let ${\cal Z}[T]$ be the subset  
$\pi_1^{-1}((\Z[T])_*^5)$ formed by those elements of $\cal Z$ whose projection to $\R^5_*$ have integer coordinates of height $\leq T$. 
By \cite{HP}, Corollary 7.2 (with no $\R^{\ell}$), there is a constant $C'_\epsilon$ such that :

\medskip

- (a') either $\pi_2({\cal Z}[T])  \subset {\cal U}_B( D \cap \mathbf \Xi[T]) \subset {\cal S}$ has less than $C'_\epsilon T^\epsilon$ elements.  Recalling the two remarks above, we then deduce  from an $o$-minimal  uniformity argument (or from a zero estimate as in  \cite{BMPZ}, Prop. 3.3)    that for some constant $C_\epsilon$, there are at most $C_\epsilon T^\epsilon$ points $\xi \in D \cap \mathbf{\Xi}$ for which $\nu U(\xi) -  \mu U^R(\xi) \in \Pi_{G_\xi}$ for some  $(\nu \neq 0, \mu) \in (\Z[T])^2$. This is Conclusion (a);

\medskip

- (b') or there is a definable connected curve ${\cal C} \subset {\cal Z}$ such that $\pi_1({\cal C}) \subset \R_*^5$   is semi-algebraic and $\pi_2({\cal C}) \subset {\cal S}$ has (real) dimension 1. Let ${\cal T} \subset D \subset X(\C)$ be the inverse image of $\pi_2({\cal C})$ under the map ${\cal U}_B$. We can view $\cal C$ as parametrized by the curve ${\cal T}$. The coordinates $\mu, \nu, \beta_0, \beta_1, \beta_2; w_0, w_1, w_2, w_0^R,  w_1^R, w_2^R$ on $\R^5 \times \R^8$, restricted to $\cal C$, then become functions of the (real) variable $\gamma \in {\cal T}$. Since $\pi_1({\cal C})$ is semi-algebraic, the functions $\mu(\gamma), \nu(\gamma), \beta_0(\gamma), \beta_1(\gamma), \beta_2(\gamma)$ generate a   field  of transcendence degree 1 (or $0$, if constant) over $\C$. In view of the incidence relations, {\it whose $\nu$-component does not vanish} by definition, the restrictions to ${\cal T}$ of the functions $w_0 = b_0, w_1 = b_1, w_2 = b_2$ generate  a field of transcendence degre $\leq 1$ over the field generated  by the restrictions to ${\cal T}$ of the functions $w^R_0 = b^R_0, w^R_1 = b^R_1, w_2 = b^R_2$. Recalling  that $U = b_0 \varpi_0 + b_1 \varpi_1 + b_2 \varpi_2$, and similarly with $U^R$, we deduce that   $U_{|{\cal T}}$ generate a field of transcendence degree $\leq 1$ over the field generated by $U^R_{|{\cal T}}$ and the $\varpi_{i|{\cal T}}$'s. By complex analyticity, the corresponding algebraic relation  extends to $D$, so $U$   generates a field of transcendence degree $\leq 1$ over the field $F^R . F_G$ generated over $\C(X)$ by $U^R$  and the $\varpi_i$'s. In view of Lemma 5, this is Conclusion (b), and the proof of Proposition 1 is completed. \qed

\section{The weakly special case over $\Z$}
From now on, we assume that {\it the sections $p$ and $q$}   are {\it linearly dependent over $End(E_0)$ modulo the subgroup  $E_0(\Q^{alg})$} of constant sections of $E_0(X)$, and look for a proof of Theorem 1.  Since its statement is invariant under multiplication of $s$ by a positive integer, and since $q$ is not constant, we can assume without loss of generality that the generic relation they satisfy takes the form 
$$p = \rho q + p_0, ~�{\rm with} �~�\rho \in End(E_0) , p_0 \in E_0(\Q^{alg}), p_0 \notin E_0^{tor}(\Q^{alg})$$
(if $p_0$ is torsion, the conclusion of Theorem 1 is trivially satisfied).  In such a case,   the initial Step 2.1.1 of the previous proof simply does not hold : contrary to the situation of Lemma 3, the set

\medskip
\centerline{ $\mathbf{\Xi}^{\Z\ell d}_{p, \delta q} = \{\xi \in X(\Q^{alg}), p(\xi)$ and $\delta q(\xi)$ are linearly dependent over $\Z$\}} 

\medskip
\noindent
may well have unbounded height. 

\medskip
  In this section, we show that if 
$$\rho = r \in \Z ,  r \neq 0,$$
bounded height for $\mathbf{\Xi}^{\Z\ell d}_{p, \delta q}$, hence for its subset  ${\mathbf{\Xi}}$, can still be recovered, thanks to    Silverman's theorem and basic orthogonality properties of N\'eron-Tate pairings. Theorem 1 then follows by reproducing most of the previous proof.

 \subsection{Bounded height}
Let again $h$ denote  the height on $X(\Q^{alg})$ attached to a divisor of degree 1.

\medskip
\noindent
{\bf Proposition 2} : {\it let $p, q \in E_0(X), p_0 \in E_0(\Q^{alg})$, $q$ not constant, and assume that there exists a} non-zero   integer  {\it $r$ such that $p= rq + p_0$. Then, there exists a constant $c'_0$ such that $h(\xi) \leq c'_0$ for any $\xi \in \mathbf{\Xi}^{\Z\ell d}_{p, \delta q}$, hence for any $\xi \in {\mathbf \Xi}$. }  

\smallskip
\noindent

\smallskip
\noindent
{\it Proof}. -  Assume for a contradiction that there exists a sequence $\xi_n, n \in \N,$ of points of $\mathbf{\Xi}^{\Z\ell d}_{p, \delta q}$ whose heights $h(\xi_n)$ tend to infinity. Denote by $\langle  \, ,   \rangle_{geo}$ the  (geometric)  N\'eron-Tate pairing on $E_0(K^{alg}) \times E_0(K^{alg})$, where $K = \Q^{alg}(X)$, and by $\langle  \, ,   \rangle_{ari}$ the  (arithmetic)  N\'eron-Tate pairing on $E_0(\Q^{alg}) \times E_0(\Q^{alg})$. 

\smallskip
Recall that for both pairings, the adjoint of $\rho \in End(E_0)$ is its complex conjugate. In particular, $\delta q(\xi) = - \overline \delta q(\xi)$ is orthogonal to $q(\xi)$, so $  \langle p(\xi_n) \, ,  q(\xi_n) \rangle_{ari} = 0$ for all $n$.  
 By Silverman (\cite{S}, or see \cite{L}, p. 306), we deduce that
$$ \langle p \, , q  \rangle_{geo} = \lim_{n \rightarrow \infty}  \frac{\langle p(\xi_n) \, ,  q(\xi_n) \rangle_{ari}}{h(\xi_n)}= 0.$$
 Now, $p = r q + p_0$, and  the constant part $E_0(\Q^{alg})$ is orthogonal to the full space $E_0(K^{alg})$ for the geometric pairing. So
$$ \langle p \, , q  \rangle_{geo} =    \langle r q \, , q  \rangle_{geo} +  \langle p_0 \, , q  \rangle_{geo} = r  \langle q \, , q  \rangle_{geo} ~�{\rm with} ~ r \neq 0.$$
Therefore, the section $q$ has vanishing N\'eron-Tate height, hence must be constant, contrary to our  hypothesis. \qed

 \subsection{Algebraic (in)dependence}
Assuming that $p = r q + p_0$ as above,  we now follow the proof of \S 2.1. All its steps go  through, except that  Conclusion (b) of Proposition 1 is now automatically satisfied. Indeed, we have $u = r v + u_0$, where $u_0 \in Lie E_0(\C)$ is a conveniently chosen elliptic logarithem of $p_0$, so $K(u)$ lies in the field $K(v) \subset F_G$, and automatically, $U = \log_G (s) $  generates a field of transcendence degree at most 1 over $F_G$. 

\medskip
To overcome this difficulty, we will now deduce from the generic relation $p = r q+ p_0$ that Conclusion (b) can here be replaced by the more precise 

\smallskip
 
- (b$^\sharp$) or the field $F_G(U)$ is algebraic  over the field $F_G(u) = F_G$

\smallskip
\noindent
(which is actually Conclusion (b2) of \S 2.1.4).

\smallskip
To check this, we use the same incidence variety ${\cal Z}$ as in \S 2.2, and follow  Alternative (b') of the discussion. Notice that  any   relation $\nu U(\xi) - \mu U^R(\xi) = \beta_0 \varpi_0(\xi) + \beta_1 \varpi_1(\xi) + \beta_2 \varpi_2(\xi)$, projected to $Lie E_0$, yields  $ \nu u(\xi) - \mu u^R(\xi) = \beta_1 \omega_1+ \beta_2 \omega_1$ 
hence since $u^R =  \delta v$ : 
$$(\nu\, r -  \mu\delta) v(\xi) = \beta_1 \omega_1 + \beta_2 \omega_2 - \nu u_0.$$
 Restricting this relation to the real curve $\cal T \subset D$,  and recalling that  $\nu \neq 0, r \neq 0$ and $\delta \notin \R$, we deduce  that if Alternative (b') holds, then the field generated over $\C$ by the restriction  of the function $v$ to   $\cal T$  lies  in  the field generated over $\C$ by the restriction  to   $\cal T$ of the real functions $\mu, \nu$ and the $\beta_i$'s, $i =  1, 2$. Since the latter field has transcendence degree at most 1 over $\C$, while $v$ is not constant, the two fields have the same algebraic closure, in which $u$ lies. The full incidence relation then implies that   $U$ is algebraic over the field $F^R.F_G(u)$ = $F_G$. This is Conclusion (b$^\sharp$). 

\smallskip

So, $\log_G( s)$ is algebraic over $F_G$. As explained in case (b2)  of \S 2.1.4,  Lemma 5.1 of \cite{BMPZ} then implies that $p$ and $q$ are linearly dependent over $End(E_0)$ and Theorem 1 is established in this ``$\rho = r \in \Z, r \neq 0$-weakly special" case.  \qed

\section{End of proof of Theorem 1}

\subsection{From weakly special to constant}
In this subsection, we assume that the projection $p \in E_0(X)$ of $s \in G(X)$ and the section $q \in E_0(X)$ are linked by a generic relation of arbitrary shape :
$$p = \rho q + p_0, ~�{\rm with} �~�\rho \in End(E_0),   p_0 \in E_0(\Q^{alg}). $$
We will deduce from the previous section that either $p$ and $q$ are linearly dependent over $End(E_0)$ (as predicted by Theorem 1), or we may assume that $\rho = 0$, i.e. $p$ itself is a constant section. 

\medskip
Replacing $s$ by $2s$  if necessary, we can write $\rho = r + r' \delta \in \Z \oplus \Z \delta \subset End(E_0)$, and consider the basic Ribet section 
$ s_{r'\alpha} = r' s^R$ of $G = G_q$ over $X$.
Its projection to $E_0(X)$ is the section $  r' p^R =    r' \delta q$. Therefore, the section $s' := s - s_{r'\alpha}$ of $G/X$ projects to 
$$\pi(s') := p' = p - r' \delta q = r q +  p_0.$$  
Moreover, for any $\xi \in X(\Q^{alg})$, $s_{r'\alpha}(\xi) = s_{r'\alpha, q(\xi)}$  is by definition a Ribet point of  $G_{q(\xi)}$.  Consequently, the set ${\mathbf \Xi} := {\mathbf \Xi}_s$ of points of  $X(\Q^{alg})$ where $s(\xi)$  is a Ribet point coincides with the set ${\mathbf \Xi}_{s'}$ similarly attached to $s'$, which is therefore infinite. Since $r \in \Z$, we deduce from the result of  \S 3 that either $p'$ and $q$, hence $p$ and $q$, are linearly dependent over $End(E_0)$, or that $r = 0$. 

\medskip
Assume now that $r = 0$, so the generic relation reads : $p = r' \delta q + p_0$, and consider again  the section $s' = s -  r' s^R$, which projects to $p' = p_0$. The corresponding set  ${\mathbf \Xi}_{s'}$ is still infinite. Therefore, we have reduced the proof of Theorem 1 to the case where $\rho = 0$, i.e. where the projection $p$ of $s$ is a constant section $p_0$. We must then show that $p_0$ is necessarily a torsion point.

\subsection{The constant case}
The word constant here refers not to the semi-abelian scheme $G/X$, which we still assume to be non constant ($q \notin E_0(\Q^{alg})$), but to the section $\pi \circ s := p = p_0 \in E_0(\Q^{alg})$. However, the duality properties of the Poincar\'e bi-extension ${\cal P}_0$ of $E_0 \times \hat E_0$ by $\G_m$  
enable us to permute the roles of $q$ and $p$, thereby   translating  the problem into one on the constant semi-abelian variety $G'_{p_0} = {\cal P}_{0 | p_0 \times \hat E_0} \in Ext(\hat E_0, \G_m)$ parametrized by the point $p_0$ of (the bidual of) $E_0$ . We must then prove that $p_0$ is torsion, i.e. that $G'_{p_0}$ is isosplit.

\medskip
Assume for a contradiction that $p_0$ is not torsion. Then for each $\xi$ in the set  $\mathbf \Xi$, there is a   relation $n p_0 - m \delta q(\xi) = 0$ with $nm \neq 0$, so $q(\xi)$ lies in the divisible hull of  $\Z. {\delta}p_0$, and is not torsion either. Consider the constant semi-abelian surface $G'_{p_0}   \in Ext (\hat E_0, \G_m)$. By  duality, we can view $s$ as a section $\check s \in G'_{p_0}(X)$, and $s(\xi)$ as a point $\check s(\xi)$ on $G'_{p_0}$ projecting to $q(\xi)$ in $\hat E_0$. Furthermore,  $\check s(\xi)$ is a 
non torsion Ribet point of $G'_{p_0}$ if and only if $s(\xi)$ is a non torsion Ribet point of $G_{q(\xi)}$~:   in the  setting of \S 1.2, this is  clear when $\varphi - \hat \varphi$ is an isomorphism, and it remains true in general via an isogeny. (In fact, it is proven in \cite{PcPbGj}, \S 5, that the  1-motive attached to $s_{\varphi, q}$ is isogenous to its Cartier dual as soon as $\varphi - \hat \varphi$ is an isogeny.)

\medskip
Therefore, the image $\check s(X)$ of $\check s$   is an irreducible curve in  $G'_{p_0}$ which contains infinitely many points of the group $\Gamma_0'$  formed  by all the Ribet points
of $G'_{p_0}$. Since this group  has finite rank (at most 1), McQuillan's Mordell-Lang theorem \cite{McQ} , as recalled in \S 1.3.(v), can be applied to  $G'_{p_0}$. We derive that $\check s$ factors through a translate by a Ribet point of  a  strict connected algebraic subgroup of $G'_{p_0}$. Since $p_0$ is not torsion, the only such one is $\G_m$, so $q(X)$ reduces to a  point of $\hat E_0$. This contradicts our assumption that $q$ is not constant, and concludes the proof of Theorem 1.  \qed
 
\subsection{Further comments}
We here  list properties of Ribet points and sections which although not used in the proof, may be relevant to further studies of unlikely intersections.

\medskip
\noindent
{\bf Remark 1} (in relation with Proposition 2) : attached to the divisor at infinity $D_\xi$ of the standard compactification of $G_{q(\xi)}$, there is a  canonical ``relative height"  $\hat h_{D_\xi}$, which vanishes on the Ribet points of $G_{q(\xi)}$ (cf. \cite{Duke}, \S 3). Is there a Zimmer-like comparison of $\hat h_{D_\xi}$ with a Weil height $h_{D_\xi}$, of the type $\hat h_{D_\xi} - h_{D_\xi} = O\big((\hat h(q(\xi)))^{1/2}\big)$,  or even just $o\big(\hat h(q(\xi))\big)$,  where $\hat h$ is the N\'eron-Tate height on $\hat E_0(\Q^{alg})$ ? Bounded height on $\mathbf \Xi$ would then follow in all cases, ``weakly special" or not. See \cite{ACL}, Thm 5.5, for an Arakelov approach to this problem.

\medskip
\noindent
{\bf Remark 2} (on the Betti maps) : let $\xi \in {\mathbf \Xi}$. By \cite{Duke}, Thm. 4, the Ribet point  $s(\xi)$  lies in the maximal compact subgroup   of  its fiber $G_{\xi}^{an}$.   So its logarithm $U(\xi)$ lies in  $ \Pi_{G_\xi} \otimes \R$, and its Betti coordinate $b_0(\xi)$ is a real number. Similarly, the Betti coordinate $b_0^R$  of the Betti map  $U^R_B$ attached to  $U^R = \log_G \circ s^R$ is actually real-valued. But a priori, not the Betti coordinate $b_0$ of $U$. It would be interesting to characterize the sections $s \in G(X)$ whose images meet the union of the maximal compact subgroups of all the fibers infinitely often. 

\medskip
\noindent 
{\bf Remark 3} (about effectivity) : as  suggested  in   \S 2.1.3.(A) (see also \S 1.3.(vii)), making the ``constants" of the text effective in terms of  the initial datas $X, E_0, q, s$,  requires a global  version of Proposition 1. One should here start with the uniformization map ${\rm Unif}�~ : \tilde {\cal P}_0 \simeq \C \rtimes (\C \times \C)\rightarrow {\cal P}_0^{an}$ of the Poincar\'e bi-extension itself, thereby  reflecting the symmetric roles played by $p$ and $q$ in the construction of Ribet sections. As far as the dependence in  $s$ is concerned, a first aim would be  to show that these constants are uniformly bounded in terms of the degree of the curve $W = s(X)$ in a projective embedding of $G$. We point out that this aim has indeed been reached in various versions of  the Mordell-Lang problem itself: see  \cite{H-P}  for a differential algebraic approach (inspired by work of Buium, and recently sharpened in \cite{Bin}), and \cite{R}, Theorem  2.4 for the general case.

\section{Appendix : Zilber-Pink for ${\cal P}_0$} 

Pink's generalization  of the  conjectures on unlikely intersections proposed by Bombieri, Masser, Zannier and by Zilber asserts :

\medskip
\noindent
{\bf Conjecture} (\cite{P}, Conjecture 1.3) : {\it let $\cal S/\C$ be a mixed Shimura variety, and let $W$ be an irreducible algebraic subvariety of $\cal S$, of dimension $d$.  Assume that the intersection of $W$ with the union of all the special subvarieties of $\cal S$ of codimension $> d$ is Zariski dense in $W$. Then, $W$ is contained in a special subvariety of $\cal S$ of positive codimension. }

\bigskip

As in the text, let again $E_0/\QQ$ be an elliptic curve with complex multiplications, with dual $\hat E_0 \simeq Ext(E_0, \G_m)$, and let ${\cal P}_0/\QQ$ be the Poincar\'e bi-extension of $E_0 \times \hat E_0$ by $\G_m$. This is a $\G_m$-torsor over $E_0 \times \hat E_0$, which admits two families of group laws. Namely, for any $q \in \hat E_0$, the restriction of ${\cal P}_0$ above $E_0 \times \{q\}$ is the semi-abelian variety   attached to $q$, viewed as a point in $Ext(E_0, \G_m)$, while  for any $p \in E_0$, the restriction of ${\cal P}_0$ above $\{p\} \times \hat E_0$ is the semi-abelian variety attached to $p$, viewed by biduality as a point in $Ext(\hat E_0, \G_m) \simeq E_0$. The important point in this Appendix is  that ${\cal P}_0$ admits a canonical structure of a mixed Shimura variety, which is described in detail in \cite{PcPbGj}. However, only a minimal knowledge of MSV theory  will be needed to prove Theorem 3 of the introduction.

\medskip
Before proving this theorem, we note (as pointed out by J. Pila) that it completely establishes Pink's conjecture for the MSV ${\cal S} = {\cal P}_0$, when the variety $W$ is defined over $\QQ$. Indeed, if the dimension $d$ of $W$ is $0$ or $3$, there is nothing to prove. If $d = 2$, then the special subvarieties of ${\cal P}_0$ of codimension $> d$ are its special points, and the statement reduces to the Andr\'e-Oort conjecture, which follows  in this case from \cite{Gao2}. So, only the case $d = 1$, i.e. Theorem 3, needs to be treated. 

\medskip
Through the first familiy of group laws above, the projection 
 $\varpi : {\cal P}_0 \rightarrow \hat E_0$ turns ${\cal P}_0$ into the universal extension ${\cal G}$ of $E_0$ by $\G_m$, over the moduli space $\hat E_0$.  For any integer $n$, we will denote by $[n]_{\cal G}$ the morphism of multiplication by $n$ of the group scheme ${\cal G}/\hat E_0$. Its  Ribet sections are well-defined, and we call their images {\it Ribet curves of ${\cal P}_0$, in the sense of ${\cal G}/\hat E_0$}. Similarly, the projection $\varpi' : {\cal P}_0 \rightarrow  E_0$ turns ${\cal P}_0$ into a group scheme ${\cal G}'/E_0$, with morphisms $[n]_{{\cal G}'}$  and {\it Ribet curves of ${\cal P}_0$, in the sense of ${\cal G}'/ E_0$}. Furthermore, $[n]_{\cal G}$ and $[n]_{{\cal G}'}$ induce the same morphism $[n]$ on the fiber $\G_m$ of $(\varpi, \varpi')$ above $(0, 0)$. With these definitions  in mind, we have the following explicit necessary conditions for an irreducible curve to be special in ${\cal P}_0$. It follows from \cite{PcPbGj}, \S 5 (see also \cite{B-E}, \S 2) that they are also sufficient, but we will not need this sharper result. 

\medskip
\noindent
{\bf Proposition 3} : {\it Let $C$ be a special curve of the MSV ${\cal P}_0$. Then,

\smallskip

i) if $\varpi : C \rightarrow \hat E_0$ is dominant, $C$ is a Ribet curve in the sense of ${\cal G}/\hat E_0$;

\smallskip

ii) if $\varpi' : C \rightarrow E_0$ is dominant, $C$ is a Ribet curve in the sense of ${\cal G}'/E_0$;

\smallskip

iii) if $(\varpi', \varpi)(C)$ is a point $(p_0, q_0)$ of $E_0 \times \hat E_0$, this point is a torsion point, and $C$ is the fiber of ${\cal P}_0$ above $(p_0, q_0)$.}

\medskip
Notice that most special curves $C$ satisfy both (i) and (ii), and are therefore Ribet curves in both senses. This reflects the self-duality of non torsion Ribet sections, already encountered in \S 4.2. As for (iii), it occurs if    neither (i) nor (ii) are satisfied.

\medskip
\noindent
{\it Proof}. - We will use the following facts, for which we refer to \cite{P}, \cite{Gao}.

\medskip
\noindent
(F1) : a point $P$ of ${\cal P}_0$ is special (if and) only if $(p,q) = (\varpi', \varpi)(P)$ is torsion in $E_0 \times \hat E_0$ and $P$ is torsion in the (isosplit) extension ${\cal G}_q$  (equivalently, in the isosplit ${\cal G}'_p$).

\smallskip

\noindent
(F2) a special curve of ${\cal P}_0$ contains a Zariski-dense set of special points, hence by F1 a Zariski-dense set of torsion points of the various fibers of ${\cal G}/\hat E_0$ (or of ${\cal G}'/E_0$).

\smallskip
\noindent
(F3) : the image of a special subvariety under a Shimura morphism (such as $(\varpi', \varpi), [n]_{\cal G}, [n]_{{\cal G}'}$) is a special subvariety.

\medskip
Let then $C \subset {\cal P}_0 = {\cal G}$ be a special curve, dominating $\hat E_0$ as in (i). By base extension 
along the finite cover  $\varpi :   X := C  \rightarrow \hat E_0$, we can view the diagonal map $X \rightarrow C_X$ as a section $s$ of the group scheme $G = {\cal G}_X := {\cal G} \times_{{\hat E}_0} X$ over $X$. We can now apply Lemma 2 of \S 1.3 (relative Manin-Mumford) to $s \in G(X)$ : by Facts F1 and F2, the set ${\bf \Xi}^{tor}_s$ is infinite and we infer that $s$ is a Ribet section of $G/X$, or factors through a torsion translate of $\G_{m/X} = \G_m \times X$. In the first case, the image $C  \subset  \cal G$ of $s(X) \subset C_X  \subset {\cal G}_X$ is a Ribet curve of ${\cal P}_0$ in the sense of ${\cal G}/\hat E_0$, as was to be shown.

\smallskip

 In the second case, a multiple $C' := [n]_{\cal G} (C)$ of $C$ lies in the fiber  $\G_m \times \hat E_0$ of ${\cal P}_0$ above $p = 0$, and is still a special curve of ${\cal P}_0$ by F3.  So, by F2, $C'$ contains  infinitely many special points   of ${\cal P}_0$ lying in $\G_m \times \hat E_0$. But by F1, these special points are contained in (in fact, fill up)  the torsion of the group $\G_m \times \hat E_0$. We can now apply the  standard Manin-Mumford theorem \cite{H} to $C' \cap (\G_m \times \hat E_0)^{tor}$,  and deduce  that $C'$ is a torsion translate of $\G_m \times \{0\}$ or of $\{1\}�\times \hat E_0$. The first conclusion cannot occur since $C'$ too dominates $\hat E_0$. So, a multiple $[m] C' = [mn]_{\cal G}(C)$ of $C$ is the image of the unit section of ${\cal G}/{\hat E_0}$. Therefore, $C$ is  in all cases a Ribet curve of ${\cal P}_0$ in the sense of ${\cal G}/\hat E_0$.
 
 \medskip
 The same proof applies to (ii), while (iii) easily follows from F1 (or from F3, in view of \cite{Gao}). This concludes the proof of Proposition 3. \qed 
 
 \bigskip
 We can now turn to the proof of Theorem 3. We will need the following complement to Fact F3.
 
 \smallskip
 \noindent
 (F4) :  under a Shimura morphism, the irreducible components of the inverse image of a special subvariety  are special subvarieties.
 
 \medskip
 So, let $W/\QQ$ be an irreducible algebraic curve in ${\cal P}_0$, which contains infinitely many points lying on  special curves   of ${\cal P}_0$. We must show that $W$ is contained in a special surface of ${\cal P}_0$. We deduce from Proposition 3 that 
 
 \smallskip
 
 (a) : $W$ contains infinitely many points lying on    Ribet curves in the sense of ${\cal G}/\hat E_0$, and if not, 
 
 (b) : ditto in the sense of ${\cal G}'/E_0$, and if not, 
 
 (c) : ditto with the fibers of ${\cal P}_0$ above the torsion points of $E_0 \times \hat E_0$.
 
 \smallskip
 
 Assume first that   $\varpi : W \rightarrow \hat E_0$ is dominant, and that we are in Case (a). Base changing along $\varpi : X = W \rightarrow \hat E_0$ as above, we may view 
 the diagonal map $X \rightarrow W_X \subset {\cal G}_X = G$ as a section $s \in G(X)$, to which Theorem 1 (or the relative Mordell-Lang Theorem 2) of \S 1.1 applies. By (a), the set ${\bf \Xi}_s$ is infinite, and we infer that the sections $p$ and $q$ attached to $s$ are linearly dependent over $End(E_0)$. So, $(\varpi', \varpi)(W)$ lies in a torsion translate of an elliptic curve $B  \subset E_0 \times \hat E_0$ passing through 0. By \cite{Gao}, these are special curves of  the MSV  $E_0 \times \hat E_0$. Therefore, by F4, $W$  lies in a special surface of ${\cal P}_0$. Vice versa, the same conclusion holds if $\varpi' : W \rightarrow E_0$ is dominant and we are in Case (b).
 
 \smallskip 
 
 Secondly, assume  that   $W$ still dominates $\hat E_0$, but that we are in Case (b). As just pointed out, we can then assume that  $W$ does not dominate $E_0$, and so, projects to a  point $p   \in E_0$ under $\varpi'$. If $p$ is not torsion,  $W$ lies in the non isosplit extension ${\cal G}'_{p} = \varpi'^{-1}(p)$ (which is then not a special surface of ${\cal P}_0)$. Now, the Ribet curves in the sense of ${\cal G}'/E_0$ meet ${\cal G}'_{p}$ at Ribet points of ${\cal G}'_{p}$, so by (b), $W$ contains infinitely many Ribet points of ${\cal G}'_{p}$. We deduce from the standard Mordell-Lang theorem  \cite{McQ} that $W$ lies in a translate of $\G_m$ by a Ribet point. But then, $W$ cannot dominate $\hat E_0$. So,  $p$ is a torsion point, and $W$ lies in $\varpi'^{-1}(p)$, which is a special surface of ${\cal P}_0$ by F4. Vice versa, the same conclusion holds if $\varpi' : W \rightarrow E_0$ is dominant and we are in Case (a).
  
 \smallskip
 Thirdly, assume that  $W$  dominates $\hat E_0$ or $E_0$, and that we are in Case (c). Then, the projection $W'$ of $W$ in $E_0 \times \hat E_0$ is a curve which contains infinitely many torsion points of $E_0 \times \hat E_0$. By Manin-Mumford, we deduce that   $W'$ lies in a torsion translate of an elliptic curve $B  \subset E_0 \times \hat E_0$ passing through 0. So, $W$ lies in a special surface of ${\cal P}_0$.
 
 \smallskip
 It remains to study the case when $W$  projects to a point $(p, q)$ of $E_0 \times \hat E_0$ under $(\varpi', \varpi)$. Then, the only special curve  of type (c) which meets $W$ is  the closure of $W$ itself, so in Case (c), $(p,q)$ is a torsion point, and $W$ lies in (many) a special surface  of  ${\cal P}_0$. Assume finally that we are in Case (a), or in Case (b). Then, $W$ contains a Ribet point of ${\cal G}_q$, or of ${\cal G}'_p$, projecting to $p \in E_0$, or to $q \in \hat E_0$. In both cases, we deduce  that the points $p$ and $q$ are linearly dependent over $End(E_0)$. So, the projection to $E_0 \times \hat E_0$ of  $W$ lies in a torsion translate of an elliptic curve $B$ passing through 0, and $W$ lies in a special surface of ${\cal P}_0$.  This concludes the proof of Theorem 3. \qed
 
\bigskip

\noindent
{\bf Acknowledgements.} -  Both authors thank the Fields Institute for invitations to their 2017 program on Unlikely Intersections, Heights, and Efficient Congruencing, where this work was initiated. It is a pleasure for the second named author to thank Gareth Jones and Jonathan Pila for motivating discussions. He would also like to thank the Engineering and Physical Sciences Research Council for support under grant  EP/N007956/1 .

 \bigskip

\noindent
{\it Authors' addresses} :

\smallskip

D.B.: daniel.bertrand@imj-prg.fr

H.S.: harry.schmidt@manchester.ac.uk

\bigskip
\noindent
{\it AMS Classification} : 14K15, 11G15, 11G50 , 11U09 

\bigskip
\noindent
{\it Key words} : semi-abelian varieties, complex multiplication, heights,  $o$-minimality, Zilber-Pink conjecture, Ribet sections.


\begin{thebibliography}{99}
  
 
  
 \bibitem{Ba} F. Barroero : CM relations in fibered powers of elliptic familes; Arxiv 1611.01955v4 and J. Inst. Math. Jussieu,  doi.org/10.1017/S1474748017000287, August 2017.
 
 \bibitem{B-C} F. Barroero, L. Capuano : Unlikely intersections in families of abelian varieties and the polynomial Pell equation; ArXiv 1801.02885v1 
 
 \bibitem{Duke} D. Bertrand : Minimal heights and polarizations on group varieties; Duke Math. J. 80 (1995),  223-250.

\bibitem{B-E} D. Bertrand : Special points and Poincar\'e bi-extensions; with an Appendix by Bas Edixhoven; ArXiv 1104.5178v1 .  

\bibitem{NDJFL} D. Bertrand : Unlikely intersections in Poincar\'e biextensions over elliptic schemes;    Notre-Dame JFL  54,  2013, 365-375.   

\bibitem{PcPbGj} D. Bertrand, B. Edixhoven :  Pink's conjecture on unlikely intersections and families of semi-abelian varieties; in preparation.


\bibitem{BMPZ}�D. Bertrand,  D. Masser, A. Pillay, U. Zannier : Relative Manin-Mumford for semi-abelian surfaces, Proc. Edin. MS 59, 2016, 837-875.

\bibitem{Bin} G. Binyamini: Bezout-type theorems for differential fields; Compos. Math. 153, 2017, 867-888.

\bibitem{ACL} A. Chambert-Loir : G\'eom\'etrie d'Arakelov et hauteurs canoniques �
sur des vari\'et\'es semi-ab\'eliennes;  Math. Ann. 314, 1989, 381-401.

\bibitem{Gao} Z. Gao : A special point problem of Andr\'e-Pink-Zannier in the universal family of abelian varieties;  Ann. Sci. SNS Pisa 17 (2017), 231-266.  

 
\bibitem{Gao2} Z. Gao : Towards the Andr\'e-Oort conjecture for mixed Shimura varieties;  J. reine   angew. Math. 732, 2017, 85-146. 

\bibitem{H} M. Hindry : Autour d'une conjecture de Serge Lang; Invent. math. 94 (1988), 575-603.

\bibitem{HP} P. Habegger, J. Pila : $o$-minimality and certain atypical intersections, Ann. Sci. ENS  49 (2016), 813-858.

\bibitem{H-P} E. Hrushovski, A. Pillay : Effective bounds for the number of transcendental points on subvarieties of semi-abelian varieties; Amer. J. Math. 122, 2000, 439-450.

\bibitem{JR} O. Jacquinot, K. Ribet : Deficient points on extensions of abelian varieties by $\G_m$; J. Number Th., 25 (1987), 133-151.

\bibitem{JS1} G. Jones, H. Schmidt : Pfaffian definitions of Weierstrass elliptic functions; ArXiv 1709.05224v3 .

\bibitem{JS2}  G. Jones, H. Schmidt : Effective relative Manin-Mumford for families of $\G_m$-extensions of an elliptic curve; in preparation.

\bibitem{JT} G. Jones, M. Thomas : Effective Pila-Wilkie for unrestricted Pfaffian surfaces;  	ArXiv:1804.08232v1.

\bibitem{L} S. Lang : Fundamental of Diophantine Geometry, Springer 1983.

\bibitem{McQ}   M. McQuillan : Division points on semi-abelian varieties, Invent. math. 120, 1995, 143-159.

\bibitem{MZ} D. Masser, U. Zannier : Torsion points on families of simple abelian varieties and Pell's equations over polynomial rings (with Appendix by V. Flynn); J. Eur. MS 17, 2012, 2375-2416.

\bibitem{P} R. Pink : A common generalization of the conjectures of Andr\'e-Oort, Manin-Mumford, and Mordell-Lang; preprint (13 p.), April 2005.

\bibitem{R} G. R\'emond: Une remarque de dynamique sur les vari\'et\'es semi-ab\' eliennes; Pacific J. Math. 254, 2011, 397-406.

\bibitem {S} J. Silverman : Heights and the specialization map for families of abelian varieties; J. reine angew. Math 342, 1983, 197-211.

\bibitem{V} E. Viada : The intersection of a curve  with algebraic subgroups in a product of elliptic curves;  Ann. Sci. SNS Pisa 11,  2003,  47-75.

\end{thebibliography}
 \end{document}